 \documentclass{ajour}
\usepackage{amsfonts,amssymb}

\begin{document}

%

\commline{Communicated by... }

 \authorrunninghead{Peter H\"ast\"o}
 \titlerunninghead{The relative metric}




\def\Rn{{\mathbb R}^n}
\def\R{\mathbb R}
\def\Rnbar{\overline{{\mathbb R}^n}}
\def\Rb{{\overline \mathbb R}}
\def\Rp{{\mathbb R}^+}
\def\C{{\mathbb C}}
\def\X{{\mathbb X}}
\def\cross{\times}
\def\ple{\preceq}
\def\pge{\succeq}
\def\ch{ {\rm \cosh} }
\def\sh{ {\rm sinh} }
\def\arsh{ {\rm arcsinh} }
\def\arch{ {\rm arccosh} }
\def\lem{{\rm {\bf Lemma.\hskip 0.5truecm}}}
\def\pro{{\bf Proposition.\hskip 0.5truecm }}
\def\cor{{\bf Corollary.\hskip 0.5truecm}}
\def\thm{{\bf Theorem.\hskip 0.5truecm}}
\def\rem{{\bf Remark.\hskip 0.5truecm}}
\def\define{{\bf Definition.\hskip 0.5truecm}}
\def\eg{{\bf Example.\hskip 0.5truecm}}
\def\lqq{\lq\lq}
\def\rqq{\rq\rq\,}
\def\rqqs{\rq\rq\hskip 0.15truecm}
\def\ineq{\not=}
\def\hbar{\vert}
\def\l{\ell}
\def\inv{ {\rm inv} } 
\def\abrv{.\hskip 0.1truecm}
\def\ident{\equiv}
\def\aeq{\approx}
\def\card{{\rm card}\,}
\def\implies{{\Rightarrow}}
\def\M{$M\hskip-0.07cm$}
\def\be{\begin{equation}}
\def\ee{\end{equation}}

\title{A New Weighted Metric: the Relative Metric I}


\author{Peter A. H\"ast\"o\thanks{Supported in part by the Academy of Finland and the Finnish Academy of Science and Letters (Viljo, Yrj\"o and Kalle V\"ais\"al\"a's Fund)}}
\affil{Department of Mathematics, P.O. Box 4 (Yliopistokatu 5), 00014, University of Helsinki, Finland}

\email{peter.hasto@helsinki.fi}





\abstract{The {\it \M--relative distance}, denoted by $\rho_M$ is a generalization of 
the $p$--relative distance, which was introduced in \cite{Li}. We establish necessary 
and sufficient conditions under which $\rho_M$ is a metric. In two special cases we derive complete characterizations of the metric. We also present a way of extending 
the results to metrics sensitive to the domain in which they are defined, thus 
finding some connections to previously studied metrics.

An auxiliary result of independent interest is an inequality related to 
Pittenger's inequality in Section 4.}

\keywords{Relative metrics, weighted metrics, Pittenger's inequality}


\begin{article}

 
\section{Introduction and main results}

In this section we introduce the problem of this study and state two useful 
corollaries of the core results. The core results themselves are stated only 
in Section 3, since they require somewhat more notation.
 The topic of this paper are 
{\it \M--relative distances}, which are functions of the form 
$$ \rho_M(x,y) := {\vert x-y\vert \over M(\vert x \vert, \vert y\vert)}, $$
where $M\colon\Rp\cross\Rp\to\Rp$ is a symmetric function satisfying 
$M(|x|,|y|)>0$ if $|x||y|>0$,
and $x$ and $y$ are in some normed space (note that $\Rp$ denotes $[0, \infty)$). 
We are interested in knowing when $\rho_M$ is a metric, 
in which case it is called the {\it \M--relative metric}. 

The first special case that we consider is when $M$ equals a power of the power mean, 
$M=A_p^q$, where
$$ A_p(x,y):=((x^p+y^p)/2)^{1/p},\ A_0(x,y):=(xy)^{1/2},$$
$$ A_{-\infty}(x,y):=\min\{x,y\} {\rm\ and\ } A_{\infty}(x,y):=\max\{x,y\},$$
for $p\in \R \setminus\{0\}$ and $x,y\in\Rp$, see also Definition \ref{ADef}. 
In this case we denote $\rho_M$ by $\rho_{p,q}$ and call it the 
{\it $(p,q)$--relative distance}. 
The $(p,1)$--relative distance was introduced by Ren--Cang Li \cite{Li}, who proved
that it is a metric in $\R$ for $p\ge 1$ and conjectured it to be one in $\C$ as well. 
Later, the $(p,1)$--relative distance was shown to be a metric in $\C$ for $p=\infty$ 
by David Day \cite{Da} and for $p\in [1,\infty)$ by Anders Barrlund \cite{Ba}.
These investigations provided the starting point for the present paper and the 
following theorem contains their results as special cases:

\begin{theorem}\label{pqMainThm} Let $q\ineq 0$. The $(p,q)$--relative distance, 
$$ \rho_{p,q}(x,y)= { \vert x-y\vert \over A_p(\vert x\vert, \vert y\vert)^q},$$
is a metric in $\Rn$ if and only if 
$0 < q \le 1$ and $p \ge \max\{ 1-q, (2-q)/3\}$.
\end{theorem}

\begin{remark} As is done in the previously cited papers, 
we define $\rho_{p,q}(0,0)=0$ even though the 
expression for $\rho_{p,q}$ equals $0/0$ in this case. 
\end{remark} 

The second special case that we study in depth is $M(x,y)=f(x)f(y)$, 
where $f\colon \Rp\to (0,\infty)$. 

\begin{theorem}\label{ffMainThm} Let $f\colon \Rp \to (0,\infty)$ 
and $M(x,y)=f(x)f(y)$. Then $\rho_M$ is a metric in $\Rn$ if and only if 
\begin{itemize}
\item[\rm{(i)}]\  $f$ is increasing,
\item[\rm{(ii)}]\  $f(x)/x$ is decreasing for $x>0$ and
\item[\rm{(iii)}]\  $f$ is convex.
\end{itemize}
\end{theorem}

(There are non-trivial functions which satisfy conditions (i)-(iii), 
 for instance the function $f(x):= (1+x^p)^{1/p}$ for $p\ge 1$.)

In the fourth section we derive an inequality of the Stolarsky mean related to 
Pittenger's inequality which is of independent interest. In the sixth section, we present a scheme for extending the results of this investigation
to metrics sensitive to the domain in which they are defined. This provides 
connections with previously studied metrics. 

This paper is the first of two papers dealing with the \M--relative distance. 
In the second paper, \cite{Ha}, we will consider various properties of the \M--relative 
metric. In particular, isometries and quasiconvexity of $\rho_M$ are studied there.


\section{Preliminaries}

\subsection{Metric and normed spaces} 

By a {\it metric} on a set $X$ we mean a function $\rho\colon X\cross X \to \Rp$
which satisfies 
\begin{itemize}
\item[1.] $\rho(x,y)$ is symmetric,
\item[2.] $\rho(x,y)\ge 0$ and $\rho(x,y)=0$ if and only if $x=y$ and
\item[3.] $\rho(x,y)\le \rho(x,z) + \rho(z,y)$ for all $x,y,z\in X$.
\end{itemize}

A function which satisfies Condition 2. is known as {\it positive definite}; 
the inequality in Condition 3. is known as the {\it triangle inequality}.

By a {\it normed space} we mean a vector space $X$ with a function $\vert \cdot\vert \colon
X \to \Rp$ which satisfies
\begin{itemize}
\item[1.] $\vert ax\vert = \vert a \vert \vert x\vert$ for $x\in X$ and $a\in \R$,
\item[2.] $\vert x \vert =0$ if and only if $x=0$ and
\item[3.] $\vert x+y\vert \le \vert x \vert + \vert y \vert$ for all $x,y\in X$.
\end{itemize}

\subsection{Ptolemaic spaces}

A metric space $(X,d)$ is called {\it Ptolemaic} if 
\be\label{Ptolemy} d(z,w) d(x,y) \le d(y,w) d(x,z) + d(x,w) d(y,z) \ee
holds for every $x,y,z,w\in X$ (for background information on Ptolemy's inequality, see e\abrv g\abrv \cite[10.9.2]{Be}). A normed space ($X$, $\vert \cdot\vert $) is Ptolemaic if 
the metric space ($X$, $d$) is Ptolemaic, where $d(x,y)=\vert x-y\vert $. The 
following lemma provides a characterization of Ptolemaic normed spaces.

\begin{lemma}\label{Pto} {\rm (\cite[6.14]{Am})} A normed space is Ptolemaic 
if and only if it is an inner product space.
\end{lemma}

Since the Ptolemaic inequality, (\ref{Ptolemy}), with $d$ equal to the Euclidean metric 
can be expressed in terms of cross-ratios (see (\ref{crDef}), in Section 6) it follows 
directly that ($\overline{\Rn}$, $q$) is a Ptolemaic metric space, where $q$ 
denotes the chordal metric,
\be\label{qDef} q(x,y):= {|x-y|\over \sqrt{1+|x|^2}\sqrt{1+|y|^2}},\ 
q(x,\infty):= {1\over \sqrt{1+|x|^2}},\ee
with $x,y\in\Rn$. The following lemma provides yet another example of a Ptolemaic space.

\begin{lemma} {\rm \cite{Ka}} Hyperbolic space is Ptolemaic. 
\end{lemma}

Thus in particular the Poincar\'e Model of the hyperbolic metric, ($B^n$, $\rho$), 
is Ptolemaic. This metric will be considered in section 5 of the sequel of
this investigation, \cite{Ha}. 

\subsection{Real-valued functions} 

An increasing function $f\colon \Rp\to\Rp$ is said to be 
{\it moderately increasing} 
(or shorter, to be MI) if $f(t)/t$ is decreasing on $(0, \infty)$. 
A function $P\colon\Rp\cross\Rp\to\Rp$ is MI 
if $P(x,\cdot)$ and $P(\cdot,x)$ are MI for every $x\in (0,\infty)$. 
Equivalently, if $P$ is symmetric and $P\not\ident 0$ then $P$ 
is MI if and only if $P(x,y)>0$ and 
$$ {z\over x} \le {P(z,y)\over P(x,y)}\le 1 \le {P(x,z)\over P(x,y)}\le 
{z\over y}$$
for all $0 < y \le z \le x$.

The next lemma shows why we have assumed that $M(x,y)>0$ for $xy>0$.

\begin{lemma} Let $P\colon\Rp\cross\Rp\to\Rp$ be symmetric and MI. 
 Then exactly one of the following conditions holds:
\begin{itemize}
\item[{\rm (i)}]\  $P\ident 0$,
\item[{\rm (ii)}]\  $P(x,y)=0$ if and only if $x=0$ or $y=0$,
\item[{\rm (iii)}]\  $P(x,y)=0$ if and only if $x=0$ and $y=0$ or
\item[{\rm (iv)}]\   $P(x,y)>0$ for every $x,y\in\Rp$.
\end{itemize}
\end{lemma}

\begin{proof} Suppose $P\not\ident 0$. Let $x,y\in(0,\infty)$ be such that $P(x,y)>0$. Then 
$$ P(z,w) \ge \min\{1, z/x\} \min\{1, w/y\} P(x,y) >0 $$ 
for every $z,w\in (0,\infty)$. Let then $x\in(0,\infty)$ be such that $P(x,0)>0$. Then 
$ P(z,0) \ge \min\{1, z/x\} P(x,0) >0 $ for every $z\in (0,\infty)$. Finally, if 
$P(0,0) >0$ then $P(x,y)>0$ for every $x,y\in\Rp$ since $P$ is increasing.
\end{proof}

\begin{lemma}\label{cont} Let $P\colon\Rp\cross\Rp\to\Rp$ be symmetric and MI. Then $P$ is continuous in $(0,\infty) \cross (0,\infty)$. 
\end{lemma}

\begin{proof} Fix points $x,y\in (0,\infty)$. 
Since $P$ is MI we have 
$$ \min\{1, z/x\}\min\{1, w/y\} P(z,w) \le P(x,y) \le $$
$$ \le \max\{1, x/z\}\max\{1, y/w\} P(z,w),$$ 
for $w,z>0$, from which it follows that $|P(x,y)-P(z,w)|$ is bounded from above by 
$$ \max\{ 1-\min\{1, z/x\}\min\{1, w/y\}, \max\{1, x/z\}\max\{1, y/w\} -1\},$$
and so the continuity is clear.
\end{proof} 

A function $P\colon \Rp\cross\Rp \to \Rp$ is said to be {\it $\alpha$--homogeneous}, $\alpha>0$, if $P(sx,sy)=s^\alpha P(x,y)$ for every $x,y,s \in\Rp$. 
A 1--homogeneous function is called {\it homogeneous}.

\begin{lemma}\label{meanLem} Let $P\colon\Rp\cross\Rp\to\Rp$ 
be symmetric, increasing and $\alpha$--homogeneous 
for some $0 < \alpha\le 1$. Then $P$ is MI.
\end{lemma}

\begin{proof} Let $x\ge z\ge y>0$. The relations 
$$xP(z,y)= x z^{\alpha}P(1,y/z) \ge z x^{\alpha}P(1,y/x)=zP(x,y)$$
and
$$yP(x,z)= y z^{\alpha}P(x/z, 1) \le z y^{\alpha}P(x/y,1)=zP(x,y)$$
imply that $P$ is MI. 
\end{proof}

\subsection{Conventions}

Recall from the introduction that $M\colon\Rp\cross\Rp\to\Rp$ 
is a symmetric function which satisfies $M(x,y)>0$ if $xy>0$. 
Throughout this paper we will 
use the short-hand notation $M(x, y):= M(\vert x \vert, \vert y\vert)$ in the 
case when $x,y\in \X$. We will denote by $\X$
a Ptolemaic normed space which is non-degenerate, i.e\abrv $\X$ non-empty and 
$\X \ineq \{0\}$. Moreover, if $M(0,0)=0$ then 
\lqq$\rho_M$ is a metric in $\X$\rqq is understood to mean that $\rho_M$ is a metric 
in $\X\setminus\{0\}$ (similarly for $\R$ or $\Rn$ in place of $\X$).


\section{The \M--relative metric}

\begin{theorem}\label{thm1} Let $M$ be MI. Then $\rho_M$ is a metric 
in $\X$ if and only if it is a metric in $\R$.
\end{theorem}

\begin{proof} Since in all cases it is clear that $\rho_M$ is symmetric and positive definite,  when we want to prove that $\rho_M$ is a metric we need to be concerned only 
with the triangle inequality. The necessity of the condition is clear; 
just restrict the metric to a one-dimensional subspace of $\X$ which 
is isometric to $\R$.

We will consider a triangle inequality of the form 
$\rho_M(x,y)\le \rho_M(x,z) + \rho_M(z,y)$.
Let $x, y, z\in \X$ be such that that $M(x,y), M(x,z), M(z,y)>0$. 
Since $M$ is increasing the case $z=0$ is trivial and we may thus assume $\vert z\vert >0$. 
For sufficiency we use the triangle inequality for the norm $\vert \cdot \vert$
and Ptolemy's inequality with $w=0$ to estimate 
$\vert x-y\vert$ in the left hand side of $\rho_M(x,y)\le \rho_M(x,z)+\rho_M(z,y)$.

We get the following two sufficient conditions for $\rho_M$ being a metric:
$$ \vert x-z\vert (1/M(x,z)-1/M(x,y))+\vert z-y\vert (1/M(z,y)-1/M(x,y))\ge 0,$$
$$ \vert x-z\vert \left( {1\over M(x,z)}- 
{\vert y \vert \over \vert z \vert M(x,y)} \right)
+ \vert z-y\vert \left({1 \over M(z,y)} - 
{\vert x\vert \over \vert z \vert M(x,y)}\right)\ge 0.$$

If $\vert z \vert \le \min \{\vert x\vert, \vert y\vert\}$, the first inequality 
holds since $M$ is increasing. The second one holds 
if $\vert z \vert \ge \max \{\vert x\vert, \vert y\vert\}$ since $f$ is MI. 
By symmetry, we may therefore assume that $\vert x\vert >\vert z\vert >\vert y \vert$. 
Then $\vert x-z\vert$ has a negative
coefficient in the first inequality, whereas $\vert z-y \vert$ has a positive one.
The roles are interchanged in the second inequality. Thus we get two 
sufficient conditions:
$$  {\vert z-y\vert \over \vert x-z\vert} \ge 
{1/M(x,y)- 1/M(x,z) \over 1/M(z,y)-1/M(x,y)}$$
and 
$$ {\vert z-y\vert \over \vert x-z\vert } \le {1/M(x,z)- \vert y \vert /
(\vert z \vert M(x,y)) \over \vert x\vert /(\vert z \vert M(x,y))- 1/M(z,y)}.$$
Now if  
$$ {1/M(x,z)- \vert y \vert /(\vert z \vert M(x,y)) \over \vert x\vert/(\vert z \vert M(x,y))- 1/M(z,y)} \ge {1/M(x,y)- 1/M(x,z) \over 1/M(z,y)-1/M(x,y)},$$
then certainly at least one of the above sufficient conditions holds. Rearranging the 
last inequality gives
\be\label{trigR} {\vert x\vert -\vert y \vert \over M(x,y)} \le 
{\vert x\vert -\vert z\vert \over M(x,z)} + 
{\vert z\vert -\vert y \vert \over M(z,y)},\ee 
the triangle inequality for $\rho_M$ in ${\R}$. Thus if $\rho_M$ is a metric in 
${\R}$, it is a metric in $\X$, so the condition is also sufficient. 
\end{proof}

\begin{remark}\label{yzx}
In the proof of Theorem \ref{thm1} we actually proved that the $\R$ in the statement 
of the theorem could be replaced by $\Rp$. Since the latter in not a vector space we 
prefer the above statement. Nevertheless, in proofs it will actually suffice to show 
that $\rho_M$ satisfies the triangle inequality for $0<y<z<x$, since the other cases 
follow from the MI condition as was seen in the proof.
\end{remark}

We may define $\rho_M$ in metric spaces as well: Let $a\in X$ be an arbitrary 
fixed point. Then we define 
$$ \rho_M(x,y) := {d(x,y) \over M(d(x,a),d(y,a))}.$$
 (As with $\X$, if $M(0,0)=0$ then we consider whether 
$\rho_M$ is a metric in $X\setminus \{a\}$.)

\begin{corollary}\label{metrCor} Let $M$ be MI and let $X$ be a 
Ptolemaic metric space and let $a\in X$ be an arbitrary fixed point. 
Then $\rho_M$ is a metric in $X$ if it is a metric in $\R$. 
\end{corollary}

\begin{proof} As in the previous proof we conclude that 
$$ {|d(x,a)-d(y,a)|\over M(d(x,a),d(y,a))} \le {|d(x,a)-d(z,a)|\over M(d(x,a),d(z,a))} + {|d(z,a)-d(y,a)|\over M(d(z,a),d(y,a))} $$
is a sufficient condition for $\rho_M$ being a metric in $X$ (this corresponds to 
(\ref{trigR})). However, since $d(x,a)$, $d(y,a)$ and $d(z,a)$ are all just real numbers 
this inequality follows from the triangle inequality of $\rho_M$ in $\R$. 
\end{proof}

\begin{corollary}\label{cor3} Let $M$ be MI. 
Then each of $\log \{1+\rho_M(x,y)\}$, $\arch \{1+\rho_M(x,y)\}$ and $\arsh \rho_M(x,y)$ 
is a metric in $\X$ if and only if it is a metric in $\R$. 
\end{corollary}

\begin{proof} Denote by $f$ the function $e^x-1$, $\cosh\{x\} - 1$ or $\sinh x$ so that the 
distance under consideration equals $f^{-1}(\rho_M)$. Applying $f$ to both sides of the triangle inequality of $f^{-1}(\rho_M)$ gives
\be\label{c1e1} \rho_M(x,y) \le \rho_M(x,z)+ \rho_M(z,y) + 
g(f^{-1}(\rho_M(x,z)), f^{-1}(\rho_M(z,y))),\ee
where $g(x,y):=f(x+y)-f(x)-f(y)$. Proceeding as in the proof of Theorem \ref{thm1}, 
we conclude that (\ref{c1e1}) follows from
\be\label{cor3E1} 
{\vert x\vert -\vert y \vert \over M(x,y)} \le {\vert x\vert -\vert z\vert \over M(x,z)} + {\vert z\vert -\vert y \vert \over M(z,y)} \ee
$$ + {\vert x\vert- \vert z\vert \over \vert x-z\vert} g\left(f^{-1}\left(\vert x-z\vert
\over M(x,z)\right), f^{-1}\left(\vert z-y\vert\over M(z,y)\right)\right).$$
We may replace the term $(\vert x\vert- \vert z\vert )/ \vert x-z\vert$ by $(\vert z\vert- \vert y\vert )/ \vert z-y\vert$ by considering the ratio $\vert x-z\vert/\vert x-z\vert$ instead of $\vert x-z\vert /\vert x-z\vert$ in the proof of Theorem \ref{thm1}. Since both conditions are sufficient we may write it as one condition by using the constant 
\be\label{logLemE1} m:=\max \left( {\vert x\vert- \vert z\vert \over \vert x-z\vert}, {\vert z\vert- \vert y\vert \over \vert z-y\vert} \right) \ge 
\sqrt{{\vert x\vert- \vert z\vert \over \vert x-z\vert} {\vert z\vert- \vert y\vert \over \vert z-y\vert}}.\ee
Then (\ref{cor3E1}) follows from the triangle inequality in $\R$ if 
$$ g\left( f^{-1}\left(\vert x\vert- \vert z\vert \over M(x,z)\right), f^{-1}\left(\vert z\vert-\vert y\vert\over M(z,y)\right)\right) \le $$
$$ \le m g\left(f^{-1}\left({\vert x-z\vert \over M(x,z)}\right), f^{-1}\left({\vert z-y\vert\over M(z,y)}\right)\right).$$ 
For $f$ equal to $e^x-1$, $\cosh\{x\} - 1$ and $\sinh x$ we see that 
$g(f^{-1}(a), f^{-1}(b))$ equals $ab$, $ab+\sqrt{a^2+2a}\sqrt{b^2+2b}$ and 
$a(\sqrt{1+b^2} -1)+ b(\sqrt{1+a^2} -1)$, respectively. Now we see that each of
these terms has either a factor $a$, $b$ or $\sqrt{ab}$, hence by choosing a suitable 
term in $m$ or the lower bound from (\ref{logLemE1}) using $\vert x-y\vert  \ge \vert x\vert -\vert y\vert $ etc\abrv the inequality follows. 
\end{proof}

The reason for considering $\log \{1+\rho_M(x,y)\}$, 
$\arch \{1+\rho_M(x,y)\}$ and $\arsh \rho_M(x,y)$ is that these metric transformations 
(see the next remark) are well-known and have been applied in various other areas, 
notably in generalizing the hyperbolic metric (see \cite[Section 5]{Ha}).
\begin{remark}\label{genMetrRem}
(i) Let $X$ be a set and $d\colon X\cross X \to \Rp$ be a function. Denote
\begin{itemize}
\item[A:] $d$ is a metric in $X$,
\item[B:] $\log \{1+ d\}$ is a metric in $X$,
\item[C:] $\arsh \{d\}$ is a metric in $X$ and
\item[D:] $\arch \{1+ d\}$ is a metric in $X$.
\end{itemize}
Then A $\implies$ B $\implies$ D and A $\implies$ C $\implies$ D, but B and C are 
not comparable, in the sense that there exists a $d$ such that B is a metric but 
C is not and the other way around. These claims are easily proved by applying 
inverse functions (that is $e^x$, $\sinh x$ and $\cosh x$) to the triangle inequality. 
For instance, to prove A $\implies$ B
we see that the triangle inequality for the $\log\{1+d\}$ variant transforms into
$1+ d(x,y) \le (1+ d(x,z))(1+d(z,y))$ which is equivalent to 
$d(x,y)\le d(x,z)+d(z,y) + d(x,z)d(z,y)$. 

(ii) Another passing remark is that if $f$ is subadditive and $d$ is a metric then 
$f\circ d$ is a metric as well. Since an MI 
function is subadditive, as noted in \cite[Remark 7.42]{AVV}, it follows, 
that the composition of an MI function with a metric is again a metric.
\end{remark}

\begin{definition} A function $P\colon \Rp\cross\Rp\to \Rp$ that satisfies 
$$ \max\{x^{\alpha}, y^{\alpha}\} \ge P(x,y) \ge \min\{x^{\alpha},y^{\alpha}\}. $$
is called an {\it $\alpha$--quasimean}, $\alpha>0$. A $1$--quasimean is called a {\it mean}. 
We define the {\it trace} of $P$ by 
$t_P(x):=P(x,1)$ for $x\in [1,\infty)$. If $P$ is an $\alpha$--homogeneous 
symmetric quasimean then 
$$ P(x,y)=y^{\alpha}P(x/y,1)= y^{\alpha} t_P(x/y) $$
for $x\ge y>0$, so that $t_P$ determines $P$ uniquely in this case.
\end{definition}

If we normalize an $\alpha$--homogeneous increasing symmetric function $P$ so that $P(1,1)=1$ then $P$ is an $\alpha$--quasimean.

\begin{definition}\label{pledef} We define a
partial order on the set of $\alpha$--quasimeans by $P\pge N$ if $t_P(x)/t_N(x)$ 
is increasing. 
\end{definition}

Note that $P\pge N$ implies that $t_P(x)\ge t_N(x)$, since 
$t_P(1)=t_N(1)=1$ by definition.

We will need the following family of quasimeans, 
related to the Stolarsky mean (see Remark \ref{stoRem}),
$$ S_p(x,y):= (1-p){ x-y \over x^{1-p}-y^{1-p}},\ S_p(x,x)=x^p,\ 0<p<1, $$
$$ S_1(x,y):=L(x,y):={x-y\over \log x - \log y},\ S_1(x,x)=x,$$
defined for $x,y\in\Rp$, $x\ineq y$. Note that $S_1(x,y)= \lim_{y\to 0} S_1(x,y)=0$ equals the classical logarithmic mean, $L$, with $S_1(x,0):=0$. 

\begin{lemma}\label{pLem} Let $0<\alpha\le 1$ and $M$ be increasing and $\alpha$--homogeneous. 
\begin{itemize}
\item[{\rm I.}]  If $M\pge S_{\alpha}$ then $\rho_M$ is a metric in $\X$.
\item[{\rm II.}]  If $\rho_M$ is a metric in $\X$, then $M(x,y)\ge S_{\alpha}(x,y)$ 
for $x,y\in\Rp$ and
$$ {M(x,1) \over S_{\alpha}(x,1)} \le {M(x^2,1) \over S_{\alpha}(x^2,1)} $$
for $x\ge 1$. 
\end{itemize}
\end{lemma}

\begin{proof} By Lemma \ref{meanLem} $M$ is MI. By Remark \ref{yzx} 
it suffices to show that the triangle inequality holds in ${\Rp}$ with $y<z<x$. 
We will consider the cases $\alpha=1$ and $\alpha<1$ separately. 

If $\alpha=1$, set $g(x):=t_M(x)/t_L(x)$ for $x\in [1,\infty)$. 
Since $M(x,0)=xM(1,0)$ and $M(z,0)=zM(1,0)$ the triangle inequality is trivial 
if $y=0$, so we may assume that $y>0$. Then the triangle inequality for $\rho_M$ becomes
\be\label{t4e1} {\log st \over g(st)}\le {\log s \over g(s)}+ {\log t \over g(t)}, \ee
where $s=x/z$ and $t=z/y$. Since $\log st = \log s + \log t$, it is clear that this 
inequality holds if $g$ is increasing, hence $L\ple M$ is a sufficient condition. 
Choosing $s=t$ shows that $g(s)\le g(s^2)$ is a necessary condition.

Assume, conversely, that $\rho_M$ is a metric. Let $0<y=x_0<x_1<... < x_{n+1}=x$ 
(note that $\X$ has a subspace isomorphic to $\R$). 
Using $M(x_i,x_{i+1}) \ge  x_i$ we conclude that
$$ {x - y \over M(x,y)} \le \sum_{i=0}^n { x_{i+1}-x_i  \over M(x_i,x_{i+1})} 
\le \sum_{i=0}^n { x_{i+1}-x_i \over x_i}, $$
and it follows by taking the limit that
$$ {x-y \over M(x,y)} \le \int_y^x {dz \over z} = \log {x\over y}$$
and hence $ L(x,y) \le M(x,y)$. 

Assume now that $\alpha<1$ and let $g(x):=t_M(x)/t_{S_{\alpha}}(x)$ 
for $x\in [1,\infty)$. If $y=0$ then the triangle inequality takes on the form 
$$ {x^{1-\alpha} - z^{1-\alpha}\over M(0,1)} \le {x-z\over M(x,z)}.$$
This is equivalent to
$$ g(x/z) \le M(1,0)/(1-\alpha) = \lim_{s\to \infty} g(s)$$
and hence follows, since $g$ is increasing. Assume then that $y>0$. 
Then the triangle inequality becomes 
$$ {x^{1-\alpha}-y^{1-\alpha} \over g(x/y)} \le {x^{1-\alpha}-z^{1-\alpha} 
\over g(x/z) }+ {z^{1-\alpha}-y^{1-\alpha} \over g(z/y) },$$
where $y<z<x$. Clearly this holds if $g$ is increasing. The necessary
conditions $g(x)\ge 1$ and $g(x)\le g(x^2)$ follow as above. 
\end{proof}

\begin{remark}\label{stoRem} For $p\in (0,1]$ and $x,y\in \Rp$ the quasimean $S_p$ 
defined above is related to Stolarsky's mean $St_{1-p}$ by
$$ St_p(x,y):= \left( {x^p-y^p \over p(x-y) }\right)^{1/(p-1)} = 
S_{1-p}(x,y)^{1/(1-p)},$$ 
for $0<p<1$ and $St_0(x,y):=L(x,y)$.
Note that the Stolarsky mean can also be defined for $p\not\in [0,1)$, 
however, we will not make use of this fact. The reader is referred to \cite{St} for 
more information on the Stolarsky mean. 
\end{remark}

\begin{remark} Strong inequalities, i.e\abrv inequalities of the type $A\pge B$, have been 
recently proved by Alzer for polygamma function \cite{Al}. Also, although not stating 
so, some people have proved strong inequalities when what they actually wanted to get 
at was just an ordinary inequality. Thus for instance Vamanamurthy and Vuorinen 
proved that $AGM \pge L$, where $AGM$ denotes the arithmetic-geometric mean, 
see \cite{VV}. Thus there are potentially many other forms which can be shown 
to be metrics by means of Lemma \ref{pLem}.
\end{remark}


\section{ Stolarsky mean inequalities }

\begin{definition}\label{ADef} 
Let $x,y\ge 0$. We define the {\it power-mean of order $p$} by
$$ A_p(x,y):=\left( {x^p+y^p \over 2} \right)^{1/p}$$
for $p\in \R\setminus\{0\}$ and additionally 
$$ A_{-\infty}(x,y)=\min\{x,y\},\ A_0(x,y):= \sqrt{xy}
{\rm\ and\ } A_{\infty}(x,y)=\max\{x,y\}. $$ 
Also observe the convention $A_p(x,0)=0$ for $p\le 0$.
\end{definition}

In order to use the results of the previous section, we need to investigate
the partial order \lqq$\ple$\rqqs from Definition \ref{pledef}. 
The next result is an improvement of 
a result of Tung--Po Lin in \cite{Lin} which states that $L\le A_p$ if and only if
$p\ge 1/3$. Lin's result is implied by Lemma \ref{lem2}, since 
\lqq $\ple$\rqqs implies \lqq$\le$\rqq.

\begin{lemma}\label{lem2} $L\ple A_p$ if and only if $p\in [1/3,\infty]$.
\end{lemma}

\begin{proof} Denote $t_{A_p}$ by $t_p$. Since $t_L, t_p \in C^1$, $L\ple A_p$ is equivalent to
\be\label{l2e1} {d \log t_L(x)\over dx} \le { \partial \log t_p(x)\over \partial x}.\ee
Since 
$$ {\partial^2 \log t_p(x) \over \partial p \partial x}= {x^{p-1} \log x 
\over (x^p+1)^2} >0,$$
(\ref{l2e1}) holds for $p\ge 1/3$ if it holds for $p=1/3$. Calculating (\ref{l2e1}) 
for $p=1/3$ gives
$$ {1\over x-1} - {1\over x \log x} \le {1\over x+x^{2/3}}.$$
Substituting $x=y^3$ and rearranging gives
$$ 3\log y \le (y^3-1)(1+1/y)/(y^2+1).$$
Note that equality holds for $y=1$. It suffices to show that the derivative of the right hand side is greater than
that of the left hand side. Differentiating and rearranging leads to
$$ y^6-3y^5+3y^4-2y^3 + 3y^2 -3y+1\ge 0,$$
which is equivalent to the tautology $(y-1)^4(y^2+y+1)\ge 0$.

Since \lqq$\ple$\rqqs implies \lqq$\le$\rqq, it follows from \cite{Lin} 
that $L\not\ple A_p$ for $p<1/3$. 
\end{proof}

The previous lemma can be generalized to the quasimean case:

\begin{lemma}\label{pqLem} For $0< q \le 1$, $A_p^q \pge S_q$ if and 
only if 
$$ p \ge \max\{1-q,(2-q)/3\}. $$
\end{lemma}

\begin{proof} The claim follows from the previous lemma for $q=1$. For $0<q<1$ we 
need to show that $g(x):=(x^p+1)^{q/p}(x^{1-q}-1)/(x-1)$ is increasing for all 
 $x\ge 1$ and $p \ge \max\{1-q,(2-q)/3\}$. This is equivalent to showing that the 
logarithmic derivative of $g$ is non-negative for $x\ge 1$, i.e\abrv that 
$g'(x)/g(x)\ge 0$. Rearranging the terms, we see that this is equivalent to 
\be\label{pqLe1} q( x^p + x^{1-q})(x-1) \le (x-x^{1-q})( x^p +1).\ee
Letting $x\to \infty$ and comparing exponents, we see that this can hold only if 
$p\ge 1-q$. The other bound on $p$ comes from $x\to 1^+$, however, only after some work.

As $x\to 1^+$ ($x$ tends to 1 from above), 
both sides of (\ref{pqLe1}) tend to 0. Their first derivatives both
tend to $2q$ and the second derivatives to $2q(p+1-q)$. Only in the third 
derivatives is there a difference, the left hand side tending to 
$3q(p(p-1)+q(1-q))$ and the right hand side to $3p(p-1)q+ 2p(1-q)q-2q(1-q^2)+p(1-q)q$.
Thus the right hand side of (\ref{pqLe1}) is greater than or equal to the left at 
$1^+$ only if $3p\ge 2-q$.

We still need to check the sufficiency of the condition on $p$. Since 
$A_p \pge A_s$ for $p\ge s$, it is enough to
check $p=\max\{1-q, (2-q)/3\}$. For $q\le 1/2$ set $q=1-p$ in (\ref{pqLe1}). This gives
$(2p-1)x^p(x-1)+x - x^{2p}\ge 0$. Since the second derivative of this
function is positive, the inequality follows easily. 

Now set $q=2-3p$ in (\ref{pqLe1}). Dividing both sides by $x^{p}$ and rearranging gives
$$ g_p(x):=(3p-1)(x-x^{2p-1})+ (2-3p)(1-x^{2p}) - x^{3p-1} + x^{1-p} \ge 0.$$
Since $g_{1/3}(x)= 1-x^{2/3} -1 + x^{2/3} = 0$ and $g_{1/2}(x)= 
(x-1)/2 + (1-x)/2 - x^{1/2}+ x^{1/2}=0$, the previous inequality 
follows if we show that
$\partial^2 g_p(x)/\partial p^2\le 0$ for every $x$. Now
$$ {\partial^2 g_p(x) \over \partial p^2}  =  12(x^{2p} - x^{2p-1})\log x $$
$$ -(4(3p-1)x^{2p-1} + 4(2-3p)x^{2p} - x^{1-p} + 9 x^{3p-1})\log^2 x, $$
and hence $\partial^2 g_p(x)/\partial p^2\le 0$ is equivalent to (we divide by $x^{2p}$) 
\be\label{pqLe3} 12 (1-1/x)\le (9 x^{p-1}-x^{1-3p}+4(2-3p) + 4(3p-1)/x)\log x. \ee
We will show that inequality holds for $p=1/3$ and $p=1/2$ and that the right hand side
is concave in $p$. Hence the inequality holds for $1/3<p<1/2$ as well.

For $p=1/3$, (\ref{pqLe3}) is equivalent to
$$ x(3x^{-2/3} +1)\log x \ge 4(x-1).$$
Since 
$$ \log x \ge 4{x-1\over x+3x^{1/3}}$$
holds for $x=1$, it suffices to show that the derivative of the left hand side
is greater than that of the right hand side:
$$ {1\over x} \ge 4{ x+3x^{1/3} - (x-1)(1+ x^{-2/3}) \over (x+3x^{1/3})^2} = 
4{2x^{1/3} + 1 + x^{-2/3}\over x^2 + 6x^{4/3} + 9x^{2/3} }.$$
We set $x=y^3$ and rearrange to get the equivalent condition
$$ y^5 - 2y^3 - 4y^2 + 9y -4 = (y-1)^3(y^2+3y+4) \ge 0$$
which obviously holds. Next let $p=1/2$ in (\ref{pqLe3}). We now 
need to show that
$$ 6(x-1) \le (x + 4x^{1/2} +1)\log x$$
holds for $x\ge 1$. This follows by the same procedure as for $p=1/3$.
We still have to show that the right hand side of (\ref{pqLe3}) is concave. 
However, after we differentiate twice with respect to $p$ all that remains is
$$ 9(x^{p-1} - x^{1-3p})\log^3 x. $$ 
Clearly this is negative for $x\ge 1$ and $p\le 1/2$. 
\end{proof}

As we noted in remark \ref{stoRem}, the Stolarsky mean was introduced in \cite{St} as a 
generalization of the logarithmic mean. 
The previous lemma may be reformulated to a result of independent interest. This 
result is related to Pittenger's inequality, which gives the exact range of values of $p$ 
for which the inequality $A_p^q \ge S_q$ holds (see \cite[p\abrv 204]{Bu}). 
Note that the bounds in Pittenger's inequality equal our bounds only for 
$q\in [0,1/2]\cup \{1\}$. For $q\in (1/2,1)$, there are $p$ such that the ratio 
$A_p^q / S_q$ is initially increasing but eventually decreases, however its values are 
never below $1$.  

\begin{corollary}\label{StoCor} Let $0\le q<1$. 
For fixed $y>0$ the ratio $A_p(x,y)/ St_q(x,y)$ 
is increasing in $x\ge y$ if and only if $p\ge \max\{ q, (1+q)/3\}$. 
In particular, $A_p(x,y)\ge St_q(x,y)$ for all $x,y\in\Rp$ for the same $p$ and $q$.
\end{corollary}

\begin{proof} The claim follows directly from Lemma \ref{pqLem} 
and the relationship between $S$ and $St$ given in Remark \ref{stoRem}. 
\end{proof}

\section{Applications}

In this section we combine the results from the previous two sections to derive
our main results as to when $\rho_M$ is a metric.



\begin{proof}[of Theorem \ref{pqMainThm}] 
Assume that the triangle inequality holds for some pair $(p,q)$ with $p > 0$. Then
$$ 2= \rho_{p,q} (-1,1) \le \rho_{p,q} (-1,0) + \rho_{p,q} (0,1) =2^{1+q/p},$$
hence $q\ge 0$. 

Suppose next that $p<0$ and $q>0$. Consider the triangle inequality 
$ \rho_{p,q} (\epsilon, 1) \le \rho_{p,q} (\epsilon, 1/2) + \rho_{p,q} (1/2, 1)$ as 
$\epsilon \to 0$. Then the left hand side tends to $\infty$ like 
$2^{-q/p} (1-\epsilon)\epsilon^{-q} $ and the right hand side like 
$2^{-q/p} (1/2-\epsilon)\epsilon^{-q}$, a contradiction, for sufficiently 
small $\epsilon$. 

Suppose then that $p,q<0$. Then $\rho_{p,q} (x,0)=0$ for every 
$x\in \X$, contrary to the assumption that $\rho_{p,q}$ is a metric. 
For $p=0$ we arrive at contradictions of the 
triangle inequality by letting $z$ tend to $0$ or $\infty$ according as 
$q$ is greater or less than $0$. 

Hence only the case $p,q>0$ remains to be considered.
When $q>1$, the triangle inequality $\rho_{p,q}(x,y)\le\rho_{p,q}(x,z)+ 
\rho_{p,q} (z,y)$ cannot hold, as we see by letting $z\to \infty$.

The non-trivial cases follow from Lemmas \ref{pLem} and \ref{pqLem}: 
if $p\ge \max\{ 1-q, 2/3-q/3\}$ $\rho_{p,q}$ is a metric by the lemmas. If $p<\max\{ 1-q, 2/3-q/3\}$, the ratio in the definition of $\pge$ is decreasing in a neighborhood 
of $1$ or $\infty$ (this is seen in the proof of Lemma \ref{pqLem}). In the first case $A_p^q(x,1)< S_q(x,1)$ in $(1,a)$ for some $a>1$, contradicting the first 
condition in Lemma \ref{pLem}. In the second case 
$A_p^q(x,1)/S_q(x,1) > A_p^q(x^2,1)/S_q(x^2,1)$ holds for
sufficiently large $x$ and $\rho_{p,q}$ is not a metric
by the second condition in Lemma \ref{pLem}. 
\end{proof}



We will now consider an application of Corollary \ref{cor3}.

\begin{lemma}\label{logLem2} Let $\lambda_M\colon \X\cross\X \to \Rp$ be defined 
by the formula
$$ \lambda_M(x,y):= \log \{ 1 + \rho_M(x,y) \}.$$
Then $\lambda_{A_p/c}$ is a metric in $\X$ if $c\ge 1$ for $p\in[0,\infty]$ and 
$c\ge 2^{-1/p}$ for $p\in [-\infty,0)$. The latter bound for the constant $c$ is sharp.
\end{lemma}

\begin{proof} By Corollary \ref{cor3}, it suffices to prove the claims in $\R$ with $y<z<x$. We start by showing that $\lambda_{A_p/c}$ is a metric for
 $c\ge \max \{1, 2^{-1/p}\}$. Since the case $y=0$ is trivial we may assume that $y>0$. Denote $f(x):=t_{A_p}(x)$. The triangle inequality for $\lambda_M$,
$$ \log \{1+\rho_M(x,y)\} \le \log \{1+\rho_M(x,z)\} +\log \{1+\rho_M(z,y)\},$$
is equivalent to 
\begin{equation}\label{logLem2e1} {st-1 \over f(st)}\le {s-1 \over f(s)}+ 
{t-1 \over f(t)}+ c{s-1 \over f(s)} {t-1 \over f(t)}\end{equation}
where $s=x/z$ and $t=z/y$. Since $st-1=(s-1)(t-1)+(s-1)+(t-1)$ and since 
$f$ is increasing and greater than $1$, the 
triangle inequality surely holds if $f(st)\ge f(s)f(t)/c$. However, this 
follows directly from Chebyshev's inequality (see \cite[p\abrv 50]{Bu}) 
for $p>0$ and is trivial
for $p=0$. For $p<0$ it follows from the inequality $(1+s^p)(1+t^p)\ge 1+(st)^p$.

We will now show that we cannot choose $c < 2^{-1/p}$ for $p<0$. 
Let $s=t$ in (\ref{logLem2e1}): $(s+1)/f(s^2) \le 2/f(s) + c(s-1)/(f(s)^2)$. 
As $s\to \infty$, $f(s) \to 2^{1/p}$, 
hence at the limit $ 2^{1/p} (s+1) \le 2^{1+1/p} + c 2^{2/p}(s-1)$ which implies that 
$c\ge 2^{-1/p}$.
\end{proof}

We now consider the second special case, $M(x,y)=f(x)f(y)$. 

\begin{lemma}\label{ffLem} Let $M(x,y)=f(x)f(y)$ and assume $f(x)>0$ for $x\ge 0$. 
Then $\rho_M$ is a metric in ${\R}$ if and only if $f$ is MI 
and convex in $\Rp$. 
\end{lemma}

\begin{proof} Assume that $\rho_M$ is a metric in $\R$. Let $y=-x$ in the triangle inequality for $-x<z<x$:
$$ {2x\over f(x)^2} \le {x-z \over f(x)f(z)} + {x+z \over f(x)f(z)}= {2x \over f(x)f(z)}.$$
Hence $f(x)\ge f(z)$, i.e\abrv $f$ is increasing. If $z>x$, we get instead 
$zf(x)\le xf(z)$, i.e\abrv $f(x)/x$ is decreasing, so that $f$ is MI. 
Let now $0\le y<z<x$. Then the triangle inequality
multiplied by $f(y)f(z)f(x)$ becomes
\be\label{ffLemE1} (x-y)f(z) \le (x-z)f(y) + (z-y)f(x).\ee
But this means that $f$ is convex, \cite[p. 61]{Bu}. (Alternatively, setting 
$z=:ay+(1-a)x$ gives the more standard form of the convexity condition, 
$f(ay+(1-a)x) \le af(y) + (1-a)f(x)$.) 

Assume then conversely that $f$ is MI 
and convex in $\Rp$. Then convexity 
gives (\ref{ffLemE1}) for $0\le y<z<x$, and dividing this inequality by 
$f(y)f(z)f(x)$, we get the triangle inequality for the same $y,z,x$. However, we know from Remark \ref{yzx} that this is a sufficient condition for 
$\rho_M$ to be a metric, provided $M$ is MI.
\end{proof}

\begin{proof}[of Theorem \ref{ffMainThm}] 
If $\rho_M$ is a metric in $\Rn$ it is trivially a metric in $\R$, 
since $\Rn$ includes a subspace isometrically 
 isomorphic to $\R$. Hence the claims regarding $f$ follow from 
Lemma \ref{ffLem}.
If $f\colon \Rp\to(0,\infty)$ is MI and convex then 
$\rho_M$ is a metric in $\R$ by Lemma \ref{ffLem} and hence in $\Rn$ by
 Theorem \ref{thm1}. 
\end{proof}

We now give an example of a relative-metric family where $M$ is not a mean.  
Note that this family includes the chordal metric, $q$, as a special case ($p=2$). 

\begin{example}\label{qCor} The distance
$$ {\vert x-y\vert \over \sqrt[p]{1+\vert x\vert^p}\sqrt[p]{1+\vert y\vert^p} }$$
is a metric in $\X$ if and only if $p\ge 1$. 
\end{example}


\section{Further developments}

In this section, we show how the approach of this paper can be extended to construct 
metrics that depend on the domain in which they are defined. The method is based on
interpreting $\rho_M$ as $\rho_{M, \Rn\setminus \{0\}}$, where $\rho_{M,G}$ is a 
distance function (defined in the next lemma) that depends both on the function 
$M$ and the domain $G$. The proof of the next lemma is similar to that 
that of \cite[Theorem 3.3]{Sei}. Note that the topological operations (closure, 
boundary etc.) are taken in the compact space $\Rnbar$.

\begin{lemma}\label{supLem} Let $G\subset \Rn$ with $ G \ineq \Rn$. 
If $M$ is continuous in $(0,\infty)\times(0,\infty)$ 
and $\rho_M$ is a metric then 
$$ \rho_{M,G}(x,y) := \sup_{a\in \partial G} {\vert x-y\vert \over M(\vert x-a\vert,
\vert y-a \vert)}$$
is a metric in $G$.
\end{lemma}

\begin{proof} Clearly only the triangle inequality needs to be considered. Fix two points
$x$ and $y$ in $G$. Since $M$ is continuous and $\partial G$ is a closed set in the compact space $\Rnbar$ there exists a point $a\in \partial G$ such that 
$\rho_{M,G}(x,y)=\rho_M(x-a,y-a)$. Since 
$$ \rho_M(x-a,y-a) \le \rho_M(x-a,z-a) + \rho_M(z-a,y-a) \le 
\rho_{M,G}(x,z) + \rho_{M,G}(z,y)$$
it follows that $\rho_{M,G}$ is a metric in $G$. 
\end{proof}

\begin{remark}\label{supRem}
Let $M(x,y):=\min\{x,y\}$. Then 
$$ \rho_{M,G}(x,y) = \sup_{a\in \partial G} {\vert x-y\vert \over 
\min\{\vert x-a\vert, \vert y-a \vert\}} =  {\vert x-y\vert \over 
\min\{d(x), d(y)\}},$$
where $d(x)=d(x,\partial G)$. We then have
$$ \log\{ 1 + \rho_{M,G}(x,y)\} = j_G(x,y):= 
\log\left( 1+ {|x-y|\over \min\{d(x), d(y)\}}\right),$$
which provides our first connection to a well-known metric ($j_G$ 
occurs in e.g\abrv \cite{AVV}, \cite{Sei} and \cite{Vu}).
\end{remark}

The previous lemma provides only a sufficient condition for $\rho_{M,G}$ being a metric.
It is more difficult to derive necessary conditions, but with some restrictions 
on $G$, such as convexity, this might not be impossible. 

If $M$ is homogeneous, we have a particularly interesting special case, 
as we may set
$$ \rho'_{M,G}(x,y) = \sup_{a,b\in \partial G} { \vert y, a, x, b\vert \over 
t_M (\vert x, b,a,y \vert)} = \sup_{a,b\in \partial G} {1 \over M( \vert x,y,a,b\vert, \vert x,y,b,a\vert)}$$
where 
\be\label{crDef} \vert a,b,c,d\vert :={q(a,c)q(b,d)\over q(a,b) q(c,d)}\ee
denotes the cross-ratio of the points $a,b,c,d\in {\Rnbar}$, $a\ineq b$, $c\ineq d$ 
and $q$ denotes the chordal metric (defined in (\ref{qDef})). With this notation we have 

\begin{lemma}\label{genGLem} Let $G\subset \Rnbar$ with $\card \partial G \ge 2$. 
If $M$ is increasing and homogeneous and $\rho_M$ is a metric in $\Rn$ then $\rho'_{M,G}$ 
is a metric in $G$.
\end{lemma}

\begin{proof} Fix points $x$ and $y$ in $G$. There are $a$ and $b$ in the compact 
set $\partial G$ (possibly $a=\infty$ or $b=\infty$)
for which the supremum in $\rho'_M(x,y)$ is attained. By the M\"obius invariance 
of the cross ratio, we may assume that $a=0$ and $b=\infty$. Then
$\rho'_{M,G}(x,y)= \rho_M(x',y')$, where $x'$ and $y'$ are the points corresponding 
to $x$ and $y$, and we may argue as in the proof of Lemma \ref{supLem}. 
\end{proof}

\begin{corollary}\label{genGCor} Let $G\subset \Rnbar$ with $\card  \partial G \ge 2$ and
let $M(x,y) = \max\{1, 2^{1/p}\} A_{-p}(x,y)$. Then 
$$ \delta_G^p(x,y) := \log \{1 + \rho'_{M,G}(x,y)\},$$
is a metric in $G$.
\end{corollary}

\begin{proof} Follows directly from Lemmas \ref{logLem2} and \ref{genGLem}. 
\end{proof}

With this notation we have $\delta_G(x,y)= \delta_G^{\infty}$,
where $\delta_G$ is Seittenranta's cross ratio metric (\cite{Sei}). Also note that 
$$ \delta_G^{p} = \sup_{a,b\in \partial G} \log \{1 + (\vert x,a,y,b\vert^p+ 
\vert x,b,y,a\vert^p)^{1/p} \}$$
actually takes on a quite simple form.

Instead of taking the supremum over the boundary we could integrate over it:
$$ {\tilde \rho}^p_{M,G}(x,y):= \left(\int_{\partial G} 
\rho_M(x-a, y-a)^p d\mu \right)^{1/p}$$
(defined for $\mu$--measurable $\partial G$). This metric takes the boundary into account 
in a more global manner, but is difficult to evaluate for most $G$'s.

\begin{lemma}\label{intLem} Let $\rho_M$, $G$ and $\mu$ be such that 
${\tilde \rho}^p_{M,G}(x,y)$ exists for all $x,y\in G$. If $\rho_M$ is a metric then 
${\tilde \rho}^p_{M,G}$ is a metric in $G$ for $p\ge 1$.
\end{lemma}

\begin{proof} From Minkowski's inequality  
$$ \left(\int_{ \partial G} (f+g)^p d\mu\right)^{1/p} 
\le \left(\int_{ \partial G} f^p d\mu \right)^{1/p} +
\left(\int_{ \partial G} g^p d\mu \right)^{1/p}, $$
where $f,g\ge 0$ and $p\ge 1$, and the basic triangle inequality (take $f=\rho_M(x-a,z-a)$ and $g=\rho_M(z-a,y-a)$ above) 
$ \rho_M(x,y)\le \rho_M(x,z)+ \rho_M(z,y)$ it follows that ${\tilde \rho}^p_{M,G}$ also satisfies the triangle inequality. 
\end{proof}

The integral form is quite difficult to evaluate in general, however, we can 
calculate the following explicit formulae. Note that $H^2$ denotes the 
upper half-plane.

\begin{lemma} For some constants $c_t$,
$$ {\tilde \rho}^{1/(1-2t)}_{A_2,H^2}(x,y) = c_t {|x-y|\over \sqrt[t]{|x-y|^2 +4 h^2}}$$
for $0< t<1/2$, where $h$ is the distance from the mid-point of the 
segment $[x,y]$ to the boundary of $H^2$, $h:=d( (x+y)/2, \partial H^2)$. Hence 
$$ {|x-y|\over \sqrt[t]{|x-y|^2 +4 h^2}}$$
is a metric in $H^2$ for $0< t<1/2$.
\end{lemma}

\begin{proof} The formula is derived directly by integration as follows 
\begin{eqnarray*} 
{\tilde \rho}^s_{A_2,H^2}(x,y) & = & c\left( \int_{ \partial H^2} 
{ dm_1(\xi) \over (|x-\xi|^2 + |y-\xi|^2)^{s/2} }\right)^{1/s} |x-y| \\
& = & c\left( \int_{-\infty}^{\infty}
{dw \over (a^2+b^2+h^2+w^2)^{s/2} }\right)^{1/s} |x-y|
\end{eqnarray*}
where $2a := x_1-y_1$ and $2b:=x_2-y_2$ and $h$ is as above ($x_i$ refers to the 
 $i^{\rm th}$ coordinate of $x$, similarly for $y$). Let us use the 
variable substitution $w= \sqrt{a^2+b^2+h^2} z$. Then we have 
\begin{eqnarray*} 
{\tilde \rho}^s_{A_2,H^2}(x,y) & = & c \left( \int_{-\infty}^{\infty}
{\sqrt{a^2+b^2+h^2} dz \over ((a^2+b^2+h^2)(1+z^2))^{s/2} }\right)^{1/s} |x-y| \\
& = & c (|x-y|^2+4h^2)^{(1/s-1)/2} |x-y| c_s 
\end{eqnarray*}
where 
$$ c_s:= \left( \int_{-\infty}^{\infty} {dz \over (1+z^2)^{s/2}} \right)^{1/s}.$$
Note that $c_s<\infty$ for $s>1$.

The last claim follows directly from Lemma \ref{intLem}, since $\rho_{A_2}$ 
is a metric, by Theorem \ref{pqMainThm}.
 \end{proof}

\begin{remark} We saw that 
$$ \iota_s(x,y):= {|x-y| \over (|x-y|^2+4h^2)^{(1-1/s)/2}}$$
is a metric for $s>1$. We then conclude that $\lim_{s\to \infty} \iota_s$ 
exists and hence that
$$ \iota_{\infty}(x,y):={2|x-y| \over \sqrt{|x-y|^2+4h^2}} $$
is a metric also. Note that this metric is a lower bound of the hyperbolic 
metric in the half-plane, as is seen by the path-length metric method 
in \cite[Section 4]{Ha}.
\end{remark}

We may define yet 
another distance by taking the supremum over two boundary points:
$$ \rho''_{M,G}(x,y) := \sup_{a, b\in \partial G} {\vert x-y\vert \over M(\vert x-a\vert,
\vert y-b \vert)}.$$
If we assume that $M$ is increasing and continuous, this amounts to taking 
\be\label{E1}\rho''_{M,G}(x,y) = {\vert x-y\vert \over M(d(x), d(y))},\ee
where $d(x):=d(x,\partial G)$. 

One could ask whether we could construct a general theory for 
$\rho''_{M,G}$--type metrics. This would be a very interesting theory, since it 
would involve metrics taking the geometry of the domain into account which would not 
include a complicated supremum. However, this cannot, in general, be done by our techniques: the following lemma has the important 
consequence that the proof technique of Lemma \ref{supLem} cannot be extended to 
metrics of the type $\rho_{M,G}''$. 
In the following two lemmas we will use the convention that $s e_1$ is denoted 
by $s$ etc. 

\begin{lemma} Let $G:={\Rn} \setminus \{-a,a\}$ ($a>0$), $n\ge 2$ and assume that
$M$ is increasing and continuous. Then $\rho''_{M,G}$ is a metric if and only if 
$M\ident c>0$.
\end{lemma}

\begin{proof} We assume that $\rho''_{M,G}$ is a metric. Consider first the points 
$-a-r$ and $a+r$ and let $y$ be on the line joining. We may choose $y$ so that $d(y)$ varies between 0 and $r$. Then, by the triangle inequality,
\begin{eqnarray*} 
{2(r+a) \over M(r,r)} & = & \rho''_{M,G} (-r-a,r+a) \\
& \le & \rho''_{M,G}(-r-a, y) + \rho''_{M,G}(y,r+a) = {2(r+a) \over M(r,d(y))}.
\end{eqnarray*}
Consider next the points $-a + r e_2$ and $a + r e_2$ and let $y$ be on the line 
joining them. We have
$$ {2a \over M(r,r)} \le  {2a \over M(r,d(y))}$$
but now $d(y)$ varies between $r$ and $\sqrt{r^2+a^2}$. 
Hence we have $M(x,y)\le M(x,x)$ for $y\in [0, \sqrt{x^2+a^2}]$.

Let us now consider the points $x_1 := -a-s+h e_2$, $y:= a-s+h e_2$ and 
$x_2:= a+t+ he_2$, for $t\ge 0$ and $s\le a$. We have 
$$ {2a+ s+t \over M(r,d(x_2))}  =  \rho''_{M,G}(x_1,x_2) \le $$
$$ \le  \rho''_{M,G}(x_1,y) + \rho''_{M,G}(y,x_2) = {2a \over M(r,r)} + 
{ s+t \over M(r,d(x_2))},$$
where $r= \sqrt{s^2+h^2}$. From this it follows that $M(r,r)\le M(r,y)$ where
$ y= \sqrt{t^2+h^2}=\sqrt{r^2 -s^2+t^2}$. Combining the upper and lower bounds, 
we conclude that $M(x,x)=M(x,y)$ for $y\in [\sqrt{b}, \sqrt{x^2+a^2}]$, 
where $b:= \max\{0, x^2-a^2\}$. From this it follows easily that $M\ident c$. 
\end{proof}

The next idea might be to build a theory of $\rho_{M,G}''$--type metrics 
for sufficiently regular, e.g\abrv convex domains only. The following lemma 
shows that this approach does not show much promise, either. (Note that 
$B^n$ denotes the unit ball.)

\begin{lemma} Let $P\colon (0,1]\times (0,1] \to (0,\infty)$ 
be symmetric, increasing and continuous. 
Then $\rho''_{P, B^n}$ is a metric if and only if $P\ident c>0$.
\end{lemma}

\begin{proof} According to (\ref{E1})
$$ \rho''_{P, B^n}(x,y)={\vert x-y\vert \over P(d(x), d(y))} = 
{\vert x-y\vert \over P(1-\vert x\vert , 1-\vert y\vert)}.$$
Consider the triangle inequality of the points $-r$, $0$ and $r$, $0< r<1$:
$$ {2r\over P(1-r,1-r)} \le {2r\over P(1,1-r)}.$$
This implies that $P(1,s)\le P(s,s)$ for $0< s\le 1$, and, since $P$ is increasing, 
$P(1,s)=P(t,s)$ for $0< s\le t\le 1$. 

It follows that there exists an increasing function $g\colon (0,1]\to (0,\infty)$
such that $P(x,y)=:g(\min\{x,y\})$. Take points $0<y<z<x\le 1$ on the $e_1$-axis.
Then the triangle inequality
$$ {x-y\over g(y)}\le {x-z\over g(z)} + {z-y\over g(y)}$$
implies that $g(z)\le g(y)$ and since $g$ is increasing by assumption it follows that 
$g$, and hence $P$, is constant. 
\end{proof}
Since the unit ball is in many respects as regular a domain as possible, we 
see that the prospects of generalizing the theory by restricting the domain 
are not good. A better approach seems to be 
to consider $\log \{1 +  \rho''_{M, G}(x,y)\}$, since we 
know from Remark \ref{genMetrRem} that this can be a metric even though 
$\rho''_{M, G}(x,y)$ is not. The metric $j_G$ is an example of such a metric. This line 
of research seems to be the most promising further extension.

\begin{acknowledgment}
I would like to thank Matti Vuorinen for suggesting the topic of this paper to me.
I would also like to thank Glen D. Anderson, Horst Alzer and Pentti J\"arvi 
and the referee for 
their comments on various versions of this manuscript.
\end{acknowledgment}


\end{article}
\end{document}